\documentclass[11pt,draft]{amsart}

\usepackage{amssymb}
\usepackage{url}
\usepackage{amscd}
\usepackage{amsthm}
\usepackage{amsmath,mathrsfs}
\usepackage{graphicx}

\theoremstyle{plain}

\theoremstyle{remark}

\numberwithin{equation}{section}

\begin{document}

\date{}

\title[ collinear triples in permutations]
{On the number of collinear triples in permutations}
\author{ Liangpan Li}

\address{Department of Mathematics,
Shanghai Jiaotong University,
Shanghai 200240,
People's Republic of
China}

\email{liliangpan@yahoo.com.cn}

\begin{abstract}

Let $\alpha:\mathbb{Z}_n\rightarrow\mathbb{Z}_n$ be a permutation
and $\Psi(\alpha)$ be the number of collinear triples modulo $n$ in
the graph of $\alpha$. Cooper and Solymosi had given by induction
the bound $\min_{\alpha}\Psi(\alpha)\geq\lceil(n-1)/4\rceil$ when
$n$ is a prime number. The main purpose of this paper is to give a
direct proof of that bound. Besides, the expected number of
collinear triples a permutation can have is also been determined.

\end{abstract}

\subjclass[2000]{51E15, 11T99}

\keywords{finite field, collinear triple}


\maketitle

\section{Introduction}


Let $\alpha:\mathbb{Z}_n\rightarrow\mathbb{Z}_n$ be a permutation
and $\Psi(\alpha)$ be the number of collinear triples modulo $n$ in
\[\Gamma(\alpha)=\{(i,\alpha(i)):i\in\mathbb{Z}_n\},\]
the graph of
$\alpha$. Cooper and Solymosi  \cite{CooperSolymosi} had given by
induction the bound
\[\min_{\alpha}\Psi(\alpha)\geq\lceil(n-1)/4\rceil\]
when $n$ is a prime number. The main purpose of this paper is to
give a direct proof of that bound. Besides, the expected number of
collinear triples a permutation can have is also been determined.

It should be noted that Cooper and Solymosi \cite{CooperSolymosi}
conjectured
\[\min_{\alpha}\Psi(\alpha)=(n-1)/2,\]
and Cooper  had improved the above lower bound estimate in a recent
preprint \cite{Hypergraphs}. As for the classical no-three-in-line
without modulo $n$ problem, we suggest the interested reader
visiting Wikipedia for a pleasant exposition.

\section{Expected number of collinear triples}
Let $n\geq3$ be a fixed prime number throughout this paper. Since
there are $n!$ permutations in total, the expected number of
collinear triples a permutation can have is
\[{\mathscr E}(n)=\frac{\sum_{\alpha}\Psi(\alpha)}{\displaystyle n!}.\]
Choose arbitrarily three different points $i,j,k$ from
$\mathbb{Z}_n$. Since $n$ is a prime number,
the possible choices of
\[\big(\alpha(i),\alpha(j),\alpha(k)\big)\]
in $\mathbb{Z}_n\times\mathbb{Z}_n\times\mathbb{Z}_n$
such that
\[(i,\alpha(i)),(j,\alpha(j)),(k,\alpha(k))\] are
collinear is $P(n,2)$. Hence there are $P(n,2)\cdot(n-3)!$
permutations which have collinear triple at points $i,j,k$.
Consequently
\[{\mathscr E}(n)=\frac{C(n,3)\cdot P(n,2)\cdot(n-3)!}{\displaystyle n!}=\frac{n(n-1)}{6}.\]

\section{Lower bound estimates}

Let $\alpha:\mathbb{Z}_n\rightarrow\mathbb{Z}_n$ be a permutation.
Note that for every pair of points in $\Gamma(\alpha)$, the slope of
that pair must be in $1,2,\ldots,n-1$. Partition the $C(n,2)$ pairs
into classes $\{S_k\}_{k=1}^{n-1}$ according to their slopes, say
for example, every pair in $S_k$ has slope $k$. The average pairs
one class can  have is
\[
\frac{C(n,2)}{n-1}=\frac{n}{2}.
\]
Let $B_k$ ($1\leq k\leq n-1$) be the integer satisfying
\[\sharp S_k=\frac{n}{2}-0.5+B_{k}.\]
Since
\[\sum_{k=1}^{n-1}\sharp S_k=C(n,2),\]
it follows that
\[\sum_{k=1}^{n-1}B_{k}=\frac{n-1}{2},\]
and consequently
\begin{equation}\label{positive parts}
\sum_{k:B_k>0}B_{k}\geq\frac{n-1}{2}.
\end{equation}

Next suppose $B_{k}$ is a positive integer. Partition
$\mathbb{Z}_n\times\mathbb{Z}_n$ into $n$ distinct parallel lines
$\{L_s\}_{s=1}^{n}$ with common slope $k$. Partition $S_k$ into
classes $\{E_{s}\}_{s=1}^n$ according to the lines the pairs lie in.
Write
\[V_s=\sharp\big( L_s\cap\Gamma(\alpha)\big).\]
It is not hard to see that
\begin{align}
\sum_{s=1}^{n}V_s&=n, \label{sum of points}\\
\sum_{s=1}^{n}C(V_s,2)&=\sum_{s=1}^{n}\sharp E_s=\sharp S_k,
\label{sum of pairs}
\end{align}
and there are in total
\begin{equation}
\sum_{s=1}^{n}C(V_s, 3)\label{sum of triples}
\end{equation}
collinear triples with slope $k$. For $i=1,2,\ldots,n$, let
\[m_{i}=\sharp\{s:V_s=i\},\]
the number of lines such that every line intersect with
$\Gamma(\alpha)$ at exactly $i$ points. With these notations we
rewrite (\ref{sum of points}) and (\ref{sum of pairs}) as
\begin{align}
\sum_{i=1}^{n}im_i&=n, \label{sum of points second version}\\
\sum_{i=2}^{n}m_iC(i,2)&=\frac{n}{2}-0.5+B_{k}. \label{sum of pairs
second version}
\end{align}
Multiplying (\ref{sum of pairs second version}) by 2 and subtracting
it from (\ref{sum of points second version}) subsequently yields
\begin{equation*}
\sum_{i=3}^{n}(2C(i,2)-i)m_i=2B_k-1+m_1.
\end{equation*}
Now it is true time for us to estimate from below the number of
collinear triples with slope $k$. According to (\ref{sum of
triples}) there are
\begin{equation}\label{sum of triples second version}
\sum_{i=3}^{n}m_iC(i,3)
\end{equation}
collinear triples with slope $k$. Since
\[\min_{3\leq i\leq n}\frac{C(i,3)}{2C(i,2)-i}
=\min_{3\leq i\leq n}\frac{i-1}{6}=\frac{1}{3},\]
it follows that
\begin{equation}\label{key point}
\sum_{i=3}^{n}m_iC(i,3)\geq\frac{1}{3}\cdot\sum_{i=3}^{n}(2C(i,2)-i)m_i
\geq\frac{2B_k-1+m_1}{3}\geq\frac{2B_k-1}{3}.
\end{equation}
Considering the LHS of (\ref{key point}) is an integer, we improve
the above estimate slightly into
\begin{equation}\label{key point 2}
\sum_{i=3}^{n}m_iC(i,3)\geq\Big\lceil\frac{2B_k-1}{3}\Big\rceil
\geq\Big\lceil\frac{B_k}{2}\Big\rceil.
\end{equation}
Combining (\ref{key point 2}) with (\ref{positive parts}) yields
\[\sum_{k:B_k>0}\Big\lceil\frac{B_k}{2}\Big\rceil
\geq\Big\lceil\sum_{k:B_k>0}\frac{B_k}{2}\Big\rceil
\geq\Big\lceil\frac{n-1}{4}\Big\rceil.\] Hence there are at least
$\big\lceil\frac{n-1}{4}\big\rceil$ collinear triples in the graph
of $\alpha$.


\section{Remarks}
Suppose $\Gamma$ is a subset of $\mathbb{Z}_n\times\mathbb{Z}_n$ and
let $\Psi(\Gamma)$ be the number of collinear triples modulo $n$ in
$\Gamma$. Cooper and Solymosi \cite{CooperSolymosi} proved that if
$\sharp\Gamma\geq n+3$, then $\Psi(\Gamma)>0.$ In fact, both Cooper
and Solymosi'  and this paper's methods can further show that
\begin{equation}\label{best estimate}\sharp\Gamma=n+2\ \Rightarrow\
\Psi(\Gamma)\geq\Big\lceil\frac{n+1}{4}\Big\rceil,\end{equation} the
interested reader can easily provide the details. We conclude this
paper with two examples to show that (\ref{best estimate}) is best
possible in two senses when $n=5$. Define two sets in
$\mathbb{Z}_5\times\mathbb{Z}_5$ by
\begin{align*}
\Gamma_1&=\big\{(0,0),(0,1),(1,2),(1,3),(2,2),(4,1)\big\},\\
\Gamma_2&=\big\{(0,0),(0,1),(1,2),(1,3),(2,2),(4,1),(2,1)\big\}.
\end{align*}
Then $\Gamma_1$ is free of collinear triples and $\Psi(\Gamma_2)=2$.


\end{document}